\documentclass[11pt]{amsart}
\usepackage{comment}
\usepackage{amsmath}

\usepackage[usenames]{color}
\usepackage{amssymb}
\usepackage{latexsym}
\usepackage{graphics} 
\usepackage{graphicx}
\usepackage{epstopdf}  
\usepackage{enumitem}
\usepackage[all]{xy}
\usepackage{verbatim} 
\DeclareSymbolFont{AMSb}{U}{msb}{m}{n}
\DeclareMathSymbol{\N}{\mathbin}{AMSb}{"4E}
\DeclareMathSymbol{\Z}{\mathbin}{AMSb}{"5A}
\DeclareMathSymbol{\R}{\mathbin}{AMSb}{"52}
\DeclareMathSymbol{\Q}{\mathbin}{AMSb}{"51}
\DeclareMathSymbol{\I}{\mathbin}{AMSb}{"49}

\newtheorem{theorem}{Theorem}[section]

\newtheorem{prop}{Proposition}[section]

\newtheorem{corollary}{Corollary}

\newtheorem{lemma}{Lemma}[section]

\newtheorem{defi}{Definition}[section]
\newtheorem{q}[theorem]{Question}

\newtheorem*{Open} {Open Question}
\newtheorem*{wrongT}{Sawollek's Theorem}

\newtheorem*{AdamsT}{Adams et al.'s Theorem}

\newtheorem*{LickThis}{Lickorish and Thistlethwaite's Theorem}
\newtheorem*{This}{Thistlethwaite's Theorem}

\newtheorem*{tait}{Tait's Conjecture}
\theoremstyle{definition} 
\newtheorem{definition}[theorem]{Definition}

\newcommand{\bi}{\begin{itemize}}
\newcommand{\ei}{\end{itemize}}
\newcommand{\be}{\begin{enumerate}}
\newcommand{\ee}{\end{enumerate}}

\begin{document}

\begin{abstract}  We prove that all $1$-vertex spatial graphs with adequate diagrams have minimal crossing number, and that spatial graph diagrams obtained by replacing vertices and edges of a planar embedded graph by minimal crossing link or spatial graph diagrams have minimal crossing number.  Finally, we give an example in answer to a question of Adams et al. about minimal crossing diagrams of rigid vertex graphs.

\end{abstract}
\title{Minimal crossing diagrams of spatial graphs}
  \date{\today}
  
   \author[E. Flapan, H.N.\ Howards]{Erica Flapan, Hugh Howards}
    \subjclass{ 57K10, 57M15, 05C10}

    \keywords{ spatial graphs, tangles, reduced and alternating, adequate, Tait's conjecture, minimal crossing number}
    
    \address{Emerita Professor, Department of Mathematics, Pomona College, Claremont, CA 91711, USA}

\address{Department of Mathematics, Wake Forest University, Winston-Salem, NC 27109, USA}

  \maketitle

\section{Tait's Conjecture and its generalizations}
The study of minimal crossing diagrams of links began in 1898 when P. G. Tait made the following conjecture, which was independently proved by Thistlethwaite \cite{Thistle}, Kauffman \cite{Kauff}, and Murasugi \cite{Mura1} in 1987.

\begin{tait}\cite{Thistle,Kauff,Mura1}  A reduced alternating diagram of a link $L$ has minimal crossing number, and if $L$ is prime then no non-alternating diagram has minimal crossing number.\end{tait}

\begin{figure}[http]
\begin{center}
\includegraphics[scale=.4]{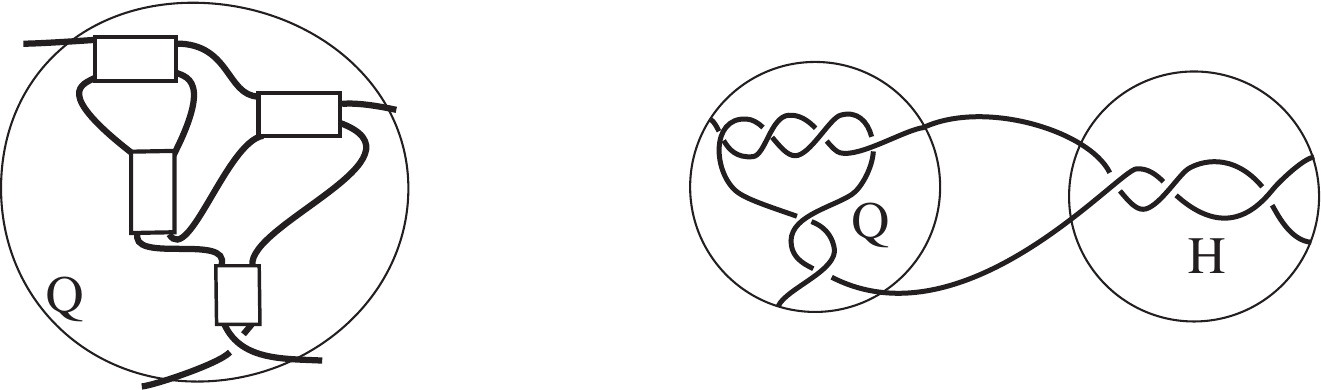}
\caption{A rational tangle can be written as $Q+H$ where $Q$ is trivial or is alternating and in the form on the left and $H$ is a (possibly trivial) row of horizontal twists.}
\label{AdequateSum}
\end{center}
\end{figure}

Lickorish and Thistletwaite \cite{LT} generalized Tait's Conjecture to Montesinos links (i.e., the numerator closure of a sum of rational tangles) with specific types of diagrams.  In particular, a rational tangle in standard alternating form can be expressed as a sum $Q+H$ where $Q$ is either trivial or an alternating rational tangle containing at least two crossings, including a final crossing between the bottom strings (as shown on the left side of Figure~\ref{AdequateSum}), and $H$ is a row of (possibly zero) horizontal twists.  Furthermore, if we have a tangle of the form $Q_1+H_1+ \dots +Q_m+H_m$, then all of the $H_i$ can be put together into a single row of horizontal twists $H$ at the right. Lickorish and Thistletwaite defined a Montesinos link diagram written as the numerator closure of a tangle $Q_1+\dots+Q_m+H$ where the $Q_i$ and $H$ are as above to be \emph{reduced} if either the diagram  is alternating or $H$ is empty.

\begin{LickThis} \cite{LT} A reduced Montesinos link diagram has minimal crossing number.
\end{LickThis}

Thistlethwaite also proved a generalization of Tait's Conjecture to adequate links, which we define below.  

\begin{figure}[h!]
 	\centering	\includegraphics[width=.4\textwidth]{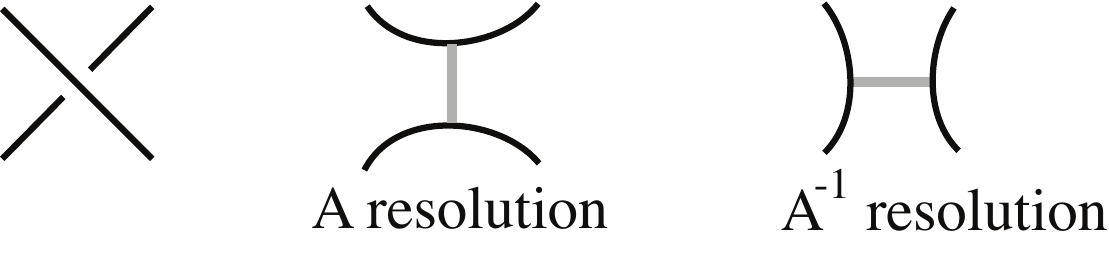}
 	\caption{$A$- and $A^{-1}$-resolutions of a crossing and their associated grey segments.}
 	\label{smoothings}
 \end{figure}

  \begin{defi} Let $G$ be a diagram of a link.  If we resolve all of the crossings as in the center of Figure~\ref{smoothings}, we obtain a collection of circles called the {\bf all-$A$ resolution} of $G$.  If we resolve all of the crossings as on the right, we obtain a collection of circles called the {\bf all-$A^{-1}$ resolution} of $G$.  \end{defi}

 \begin{defi} We add grey segments to the collection of circles in the all-$A$ or all-$A^{-1}$ resolution of a link diagram $G$ in place of crossings as in Figure~\ref{smoothings}.  We say $G$ is {\bf $A$-adequate} or {\bf $A^{-1}$-adequate} if every grey segment has its endpoints on distinct circles of the all-$A$ resolution or all-$A^{-1}$ resolution, respectively.  If $G$ is both $A$-adequate and $A^{-1}$-adequate, then we say $G$ is {\bf adequate}.
 \end{defi}

\begin{This}\cite{This} An adequate diagram of a link has minimal crossing number.
\end{This}

In addition to the above results, Tait's Conjecture has been generalized to torus links by Murasugi \cite{Mura}, to cable knots by Kalfagianni and McConkey \cite{KM}, to composite knots by Lackenby \cite{Lack}, as well as to other families of links (see for example \cite{D,HY,Stoimenow1}).

 Sawollek generalized Tait's Conjecture to reduced alternating $4$-valent spatial graphs using the definition below.

\begin{defi} \label{RedAltSawollek} \cite{Sawollek2}  Let $G$ be a diagram of a $4$-valent spatial graph.    We say $G$ is {\bf reduced and alternating}, if every smoothing of the vertices of $G$ in the plane of diagram yields a reduced alternating diagram of a link.   \end{defi}

 \begin{wrongT}\cite{Sawollek2} \label{SawollekT} Let $G$ be a reduced alternating diagram of a $4$-valent spatial graph $\Gamma$.  Then $G$ has minimal crossing number.
 \end{wrongT}

In the next section, we generalize Sawollek's result to $1$-vertex graphs with an arbitrary number of edges where \emph{ reduced alternating} is replaced by the weaker condition of \emph{adequate}.  Then, in Section~\ref{Montesinos} we generalize it to spatial graph diagrams made up of pieces with minimal crossing number lying on a planar framework.  Finally, in Section~\ref{rigid}, we present Adams et al.'s generalization of Tait's Conjecture to rigid vertex spatial graphs \cite{Adams} and give an example in answer to a question they raise.

\section{Adequate $1$-vertex graphs}\label{ourTait}

Below we extend the definitions of reduced, alternating, and adequate from links to spatial graphs with arbitrary even valence.  

\begin{defi} \label{2nReducedAlternating} Let $G$ be a diagram of an even valent spatial graph.    A {\bf smoothing} of $G$ is a link diagram obtained by replacing neighborhoods of vertices by non-intersecting planar segments joining the endpoints of edges.  \end{defi}

\begin{defi} Let $G$ be a diagram of an even valent spatial graph.  We say $G$ is {\bf reduced}, {\bf alternating}, or {\bf adequate} if every smoothing of the vertices of $G$ yields a link diagram which is reduced, alternating, or adequate, respectively.
\end{defi}

\begin{figure}[h!]
\includegraphics[scale=.5]{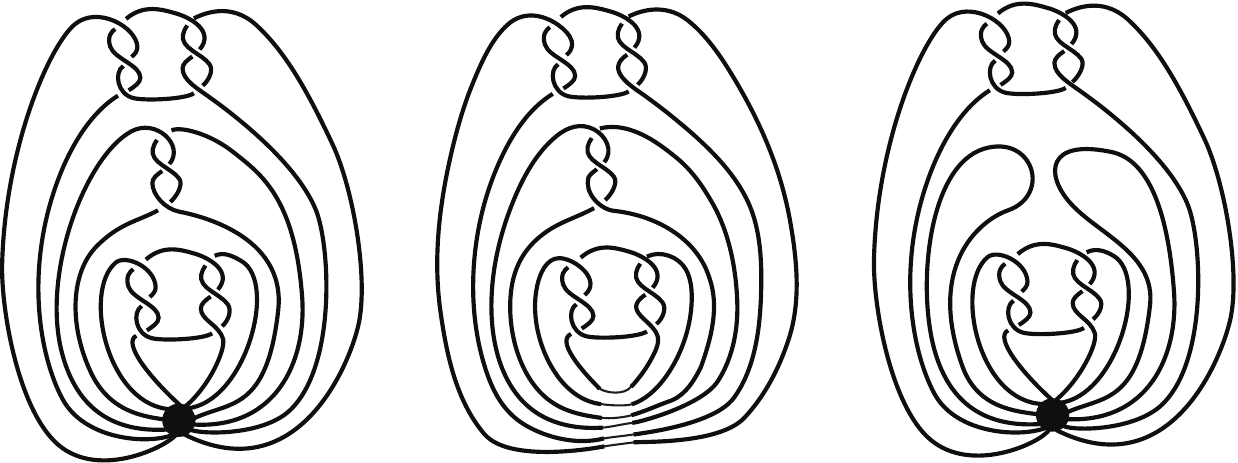}
\caption{The smoothing in the center image is not reduced.  On the right is a diagram with fewer crossings than the original.}
\label{badsmoothing}
\end{figure}

Observe that the diagram of the spatial graph on the left in Figure~\ref{badsmoothing} has the property that connecting adjacent endpoints would give us a reduced alternating link.  However, the graph is not reduced because there is a different smoothing of the vertex (illustrated in the center) which is not reduced.  Furthermore, the diagram on the left in Figure~\ref{badsmoothing} does not have minimal crossing number, since we can twist around the vertex to get a diagram with fewer crossings, illustrated on the right.

\begin{figure}[h!]
\includegraphics[scale=.5]{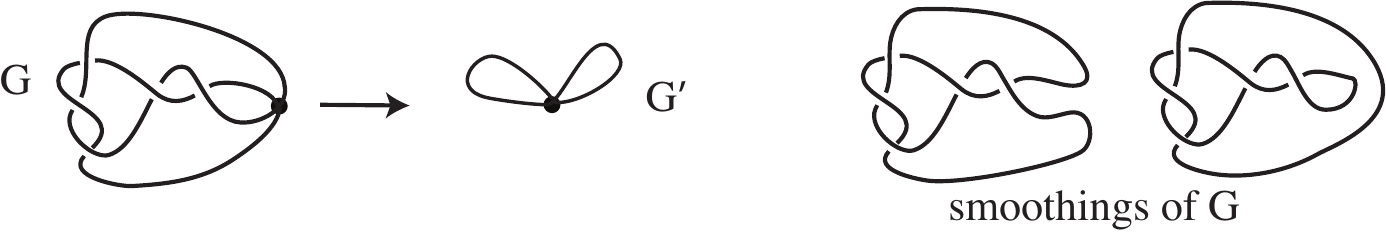}
\caption{$G$ is isotopic to $G'$, but neither smoothing of $G$ is isotopic to a smoothing of $G'$.}
\label{smoothing}
\end{figure}

It is tempting to think that if a spatial graph diagram $G$ can be isotoped to a diagram $G'$, then a link obtained by a smoothing of the vertices of $G$ will be isotopic to a link obtained by a smoothing of the vertices of $G'$. However, the spatial graph diagram $G$ on the left in Figure~\ref{smoothing} is isotopic $G'$, but neither smoothing of $G$ is isotopic to a smoothing of $G'$.   In particular, both smoothings of $G$ are non-trivial knots, but neither smoothing of $G'$ is a non-trivial knot.

Thus instead of using smoothings of a spatial graph diagram to obtain a link we use $n$-string tangles to associate a different link to the spatial graph.

\begin{definition} An {\bf $n$-string tangle diagram} is a diagram of a disk $B$ together with a set $T$ of $n$ arcs where the endpoints of the arcs on $\partial B$ are numbered $1$, \dots, $2n$ going clockwise.  Such a tangle diagram is said to be {\bf adequate} if no matter how we join the endpoints of $T$ by disjoint planar arcs we get an adequate link diagram.   \end{definition}

\begin{definition} Let $(B,T)$ be an $n$-string tangle diagram, and let $(B^*,T^*)$ be the mirror image of $(B,T)$ so that the endpoints of $T^*$ are numbered $1$, \dots, $2n$ going counterclockwise.  Let $L_T$ be the link diagram obtained by gluing $B$ and $B^*$ together so that each endpoint $j$ of $T$ is attached to the endpoint $j$ of $T^*$. 
\end{definition}

\begin{lemma} \label{Lem:LG}  Let $(B,T)$ be an adequate $n$-string tangle diagram.  Then the link diagram $L_T$ is adequate. \end{lemma}

\begin{proof} Suppose without loss of generality, that the all-$A$ resolution of $L_T$ has a grey segment $S\subseteq B$ with both endpoints on a single circle $C$.

Let $R$ denote the set of arcs of the all-$A$ resolution of $L_T$ which are contained in $B^*$ and have their endpoints on $\partial B$. Let $J$ be the link diagram $T\cup R$.  The all-$A$ resolution of $J$ intersected with $B$ is identical to the all-$A$ resolution of $L$ intersected with $B$.  Since the grey segment $S\subseteq B$, either the circle $C$ is contained in $B$ or $C$ is the union of arcs in $B$ and arcs in $R$.  In either case, $C$ is a circle in the all-$A$ resolution of $J$.  

But since $T$ is an adequate tangle and $R$ is a set of disjoint planar arcs that join the endpoints of $T$, by definition $J$ must be adequate.  As this is impossible, it follows that $L_T$ is adequate.   \end{proof}

\begin{defi} \label{tangle} Let $G$ be a spatial graph diagram $G$ with a single vertex $v$.  Define the {\bf  tangle diagram associated to $G$} to be $(B,T)$ where $B=\mathrm{cl}(S^2-N(v))$ and $T=G\cap B$.  \end{defi}

\begin{lemma}  \label{MirrorTangle} Let $G$ be a spatial graph diagram with a single vertex and suppose $G'$ is isotopic to $G$.  Let $T$ and $T'$ be the tangle diagrams associated to $G$ and $G'$, respectively.  Then the link diagrams $L_T$ and $L_{T'}$ are isotopic.
\end{lemma}

\begin{figure}[h!]
\includegraphics[scale=.28]{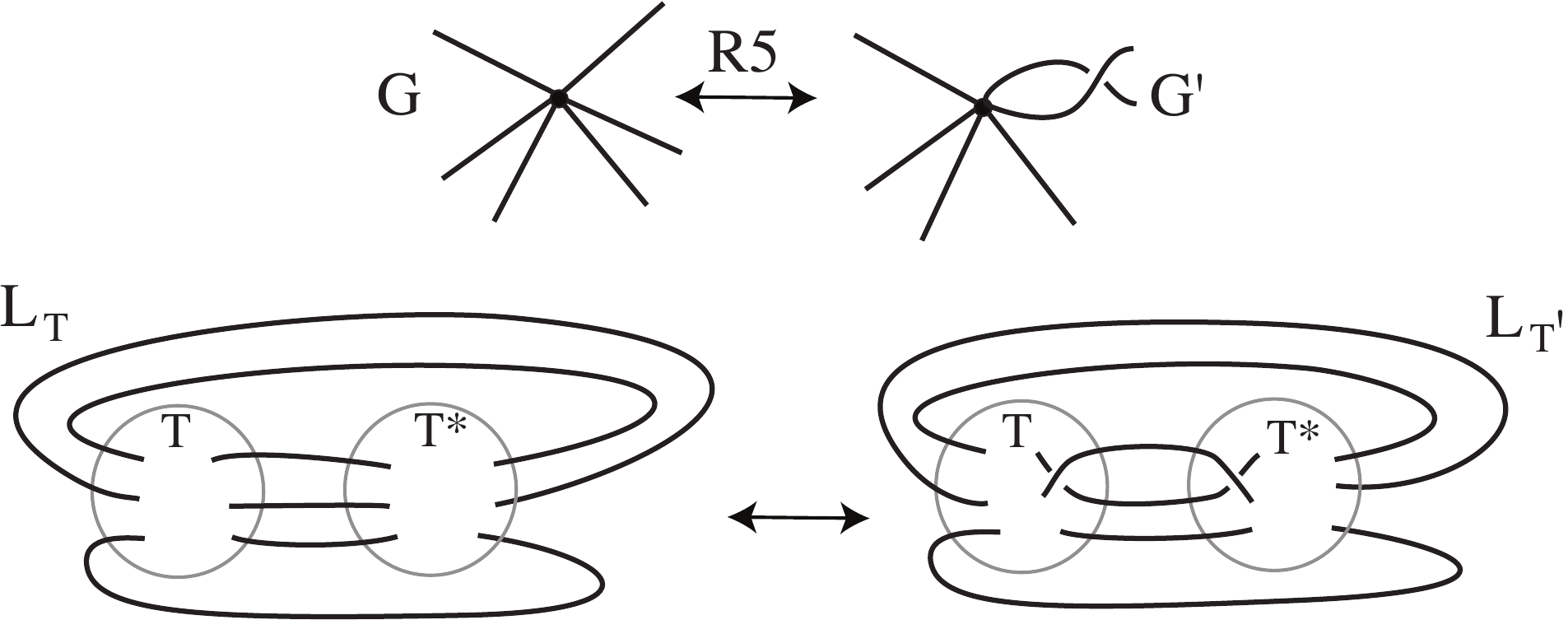}
\caption{An R5 move taking $G$ to $G'$ induces an isotopy taking $L_T$ to $L_{T'}$.}
\label{twisting}
\end{figure}

\begin{proof}   Suppose that one of the five Reidemeister moves for spatial graphs takes $G$ to $G'$.  We see below that there is an isotopy taking $L_T$ to $L_{T'}$.  

Let $B=\mathrm{cl}(S^2-N(v))$, $T=G\cap B$, and $(B^*,T^*)$ be the mirror image of $(B,T)$.  A planar isotopy or an  R1, R2, or R3 move taking $G$ to $G'$ induces an analogous move on $(B,T)$.  On $L_G$, this does an R1, R2, or R3 move on $T$ and $T^*$, giving us an isotopy from $L_T$ to $L_{T'}$.

An R4 move taking $G$ to $G'$ pulls an arc over or under the vertex $v$.  This means we are pulling an arc of $T$ over or under $N(v)$. On $L_T$, this pulls an arc of $T$ along the inside of $\partial B$ in $(B,T)$ and along the inside of $\partial B^*$ in $(B^*,T^*)$.  (In Figure~\ref{R4}, the arc being pulled is  turquoise.) This gives us an isotopy from $L_T$ to $L_{T'}$.  

\begin{figure}[h!]
\includegraphics[scale=.28]{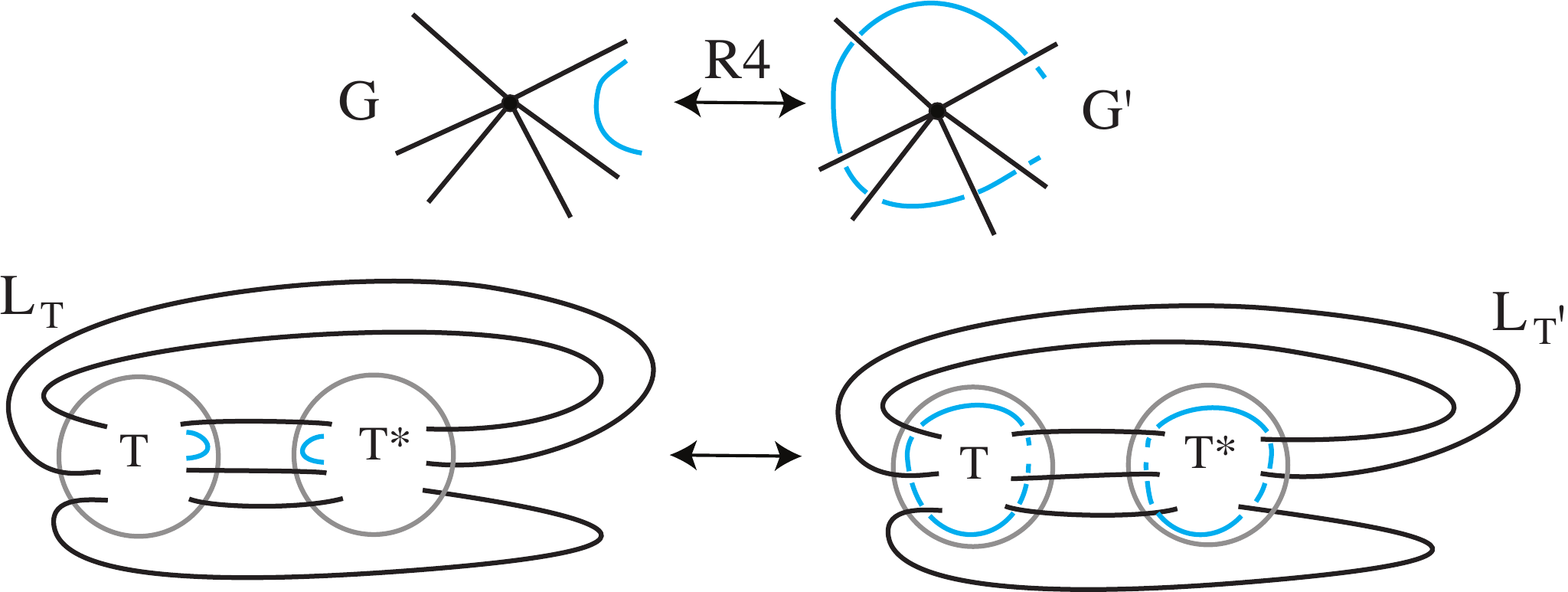}
\caption{An R4 move taking $G$ to $G'$ induces an isotopy taking $L_T$ to $L_{T'}$.}
\label{R4}
\end{figure}

Finally, an R5 move taking $G$ to $G'$, twists or untwist a pair of adjacent arcs around $v$.  On $L_T$, this adds or removes a crossing to $(B,T)$ near $\partial B$ and adds or removes its mirror image to $(B^*,T^*)$ near $\partial B^*$.  As we see in Figure~\ref{twisting}, this gives us an R2 move taking $L_T$ to $L_{T'}$.  \end{proof}

\begin{lemma} \label{lemma:onevertex} Let $G$ be a spatial graph diagram with a single vertex with associated tangle diagram $T$. Suppose that $L_T$ is adequate.    Then $G$  has minimal crossing number.
\end{lemma}

\begin{proof}Suppose there is an isotopy taking $G$ to $G'$, and $T'$ is the associated tangle diagram of $G'$. Then by Lemma~\ref{MirrorTangle}, there is an isotopy taking $L_T$ to $L_{T'}$.  But $L_T$ and $L_{T'}$ have two times the number of crossings of $G$ and $G'$, respectively.  Since $L_T$ is adequate, Thistlethwaite's Theorem \cite{This} implies that it has minimal crossing number.  Thus $G'$ cannot have fewer crossings than $G$.  \end{proof}

\begin{theorem} \label{theorem:onevertex} Let $G$ be an adequate spatial graph diagram with a single vertex.    Then $G$  has minimal crossing number.
\end{theorem}

\begin{proof} Since $G$ is adequate, the tangle diagram $T$ associated to $G$ is adequate.   Now it follows from Lemma~\ref{Lem:LG} that the link diagram $L_T$ is adequate.  Thus the result follows by Lemma~\ref{lemma:onevertex}.\end{proof}

Observe that a reduced alternating graph is adequate, and hence Theorem~\ref{theorem:onevertex} generalizes Sawollek's Theorem to $1$-vertex graphs with any valence.

\begin{prop} \label{uncrossededge} Let $G$ be a spatial graph diagram with a single vertex and an edge $E$ with no crossings such that $G-E$ is adequate. Then $G$ is adequate, and hence achieves its minimal crossing number.
\end{prop}

\begin{proof}  Suppose $G$ is not adequate.  Then there is a smoothing of the vertex $v$ which produces a non-adequate link diagram $L$.  Without loss of generality, the all-$A$ resolution of $L$ has a grey segment $S$ with both endpoints on a single circle of the resolution.   Let $R$ be the set of arcs with endpoints in $\partial B$ obtained by intersecting the all-$A$ resolution of $L$ with $B$.

 \begin{figure}[h!]
 	\centering	\includegraphics[width=.98\textwidth]{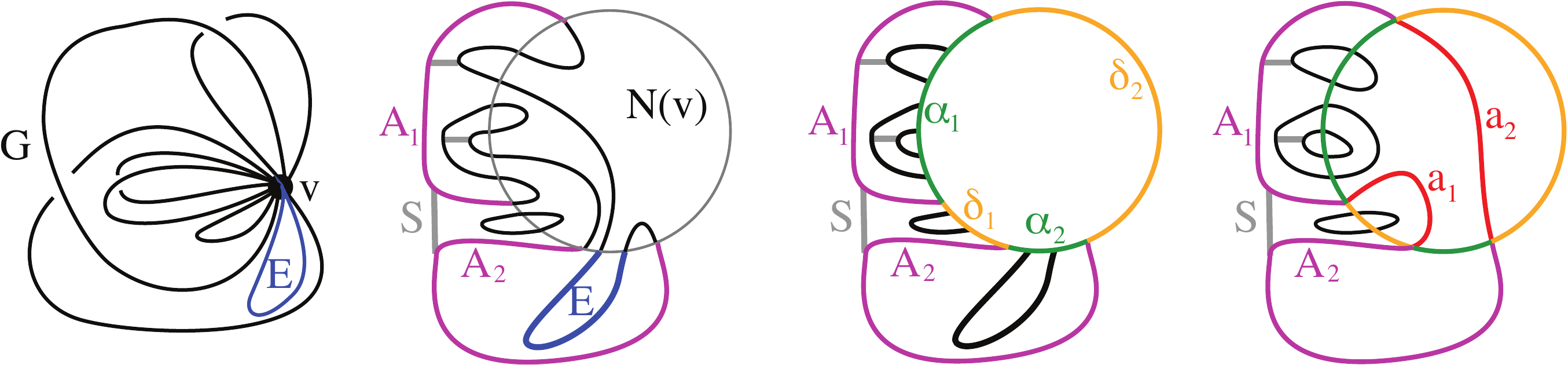}
 	\caption{ The second picture is the all-$A$ resolution of $L$, and the fourth picture is a smoothing of $v$ for $G-E$ which joins the endpoints of $A_1$ and $A_2$.}
 	\label{ExtraEdge}
 \end{figure}

Since $E$ has no crossings, $S$ must come from the resolution of a crossing in $G-E$. If $S$ has both endpoints on a single arc in $R$ or circle of the all-$A$ resolution contained in $B=\mathrm{cl}(S^2-N(v))$, then $S$ would have both endpoints on a circle in the all-$A$ resolution of any link obtained by smoothing $v$ for $G-E$.   This is impossible, since $G-E$ is adequate.  Thus $S$ must go between distinct arcs $A_1$ and $A_2$ in $R$ (see $A_1$, $A_2$, and $S$  in the second picture of Figure~\ref{ExtraEdge}).



Now since $A_1$ and $A_2$ are disjoint, there are disjoint arcs $\alpha_1$ and $\alpha_2$ in $\partial B$ such that both endpoints of $\alpha_i$ match those of $A_i$ (see the green arcs in the third picture of Figure~\ref{ExtraEdge}).  Let $\delta_1$ and $\delta_2$ denote the arcs of $\partial B-(\alpha_1\cup \alpha_2)$.   Since $S\cup A_1\cup A_2$ separates $B$ and is disjoint from the other arcs of $R$, both endpoints of each arc of $R-(A_1\cup A_2)$ are in the interior of a single $\alpha_i$ or $\delta_i$.  (In the third picture of Figure~\ref{ExtraEdge}, the arcs of $R-(A_1\cup A_2)$ are black.)

Now draw arcs $a_1$ and $a_2$ in $N(v)$ parallel to $\delta_1$ and $\delta_2$ respectively, which join the endpoints of $A_1$ with those of $A_2$. (In the fourth picture of Figure~\ref{ExtraEdge} $a_1$ and $a_2$ are red.)  Then choose disjoint arcs in $N(v)-(a_1\cup a_2)$ joining the endpoints of each arc in $R-(E\cup A_1\cup A_2)$ into a circle. This set of arcs together with $a_1$ and $a_2$ gives us a smoothing of $v$ for $G-E$, which produces a link whose all-$A$ resolution includes the circle $A_1\cup a_1\cup A_2\cup a_2$.  As $S$ has both endpoints on this circle, this contradicts our assumption that $G-E$ is adequate.  \end{proof}

\medskip

\begin{figure}[h!]
\includegraphics[scale=.7]{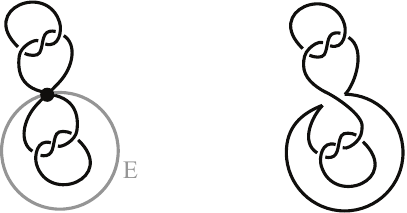}
\caption{This diagram is not alternating, but if we delete $E$ it becomes alternating.}
\label{nonalt}
\end{figure}

  Observe that the conclusion of Proposition~\ref{uncrossededge} is false if we replace adequate by reduced and alternating.  In particular, the diagram on the left in Figure~\ref{nonalt} is not alternating, as we see from the smoothing on the right.  But if we delete the edge $E$ we obtain a reduced alternating diagram.

\section{Spatial graph diagrams on planar frameworks}\label{Montesinos}

In this section, we prove that spatial graph diagrams obtained by replacing vertices and edges of a planar embedded graph by minimal crossing link or spatial graph diagrams have minimal crossing number.

\begin{defi}  Let $P$ be a planar embedded graph (i.e., a diagram with no crossings).  
\begin{itemize}
\item An {\bf edge replacement} of $P$ is obtained by replacing an edge $E_i$ by a link or spatial graph diagram $G_i$ containing the vertices of $E_i$ such that $G_i$ lies in a  neighborhood $N(E_i)$.  
\item A  {\bf vertex replacement} of $P$ is obtained by replacing a vertex $V_j$ together with its incident edges $E_{j1}$, \dots, $ E_{jn}$ by a link or spatial graph diagram $G_j$ containing each vertex $V_i$  of $E_{ji}-V_j$ such that $G_j$ lies in a neighborhood $N(V_j\cup E_{j1}\cup \dots \cup E_{jn})$.  

\item If $G$ is obtained from $P$ by some number of edge and vertex replacements, we say that  {\bf $G$ has framework $P$}.
\end{itemize}
\end{defi}

\begin{figure}[h!]
 	\centering	\includegraphics[width=.8\textwidth]{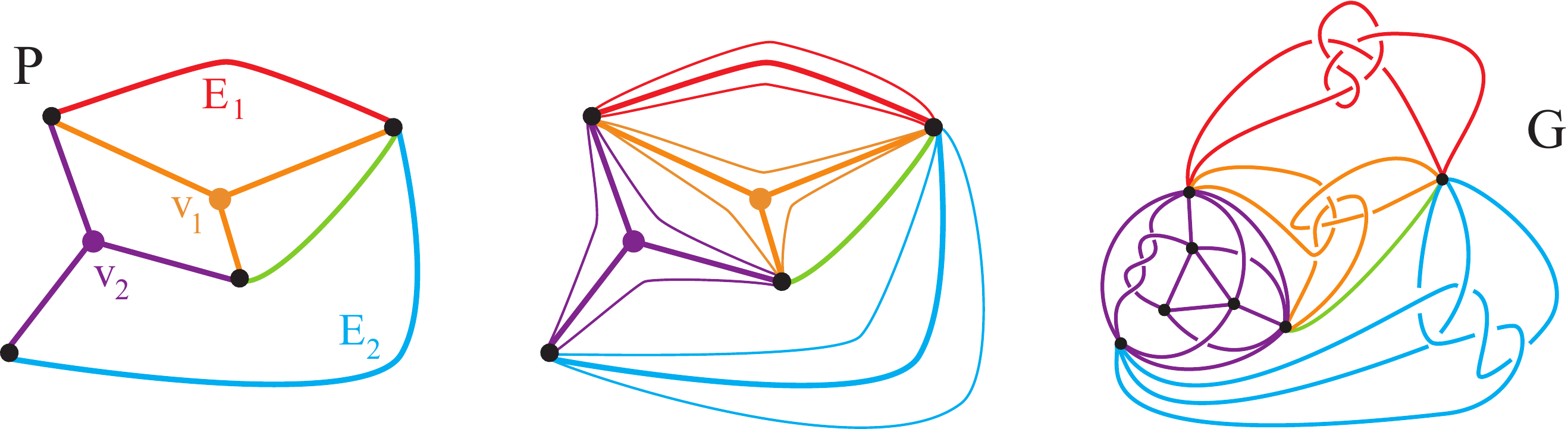}
 	\caption{On the left we illustrate $P$, in the center we show the neighborhoods of edges and vertices that we will replace, and on the right we have $G$.}
 	\label{Framework}
 \end{figure}

Figure~\ref{Framework} shows an example where edges $E_1$ and $E_2$ and vertices $v_1$ and $v_2$ are replaced by link or spatial graph diagrams of the corresponding color.  The center image shows the neighborhoods of the edges and vertices that are being replaced.  Observe that $G$ does not have even valence.

\begin{theorem} \label{theorem:framework} Let $G$ be a spatial graph diagram with a planar graph framework $P$, such that all of the link or spatial graph diagrams replacing edges and vertices have minimal crossing number.  Then $G$ has minimal crossing number.
\end{theorem}

\begin{proof} Let $G_1$, \dots, $G_m$ denote the link or spatial graph diagrams that replace edges and vertices of $P$.  Then each $G_i$ has minimal crossing number.  Since $G$ is on a planar graph framework and $G_1$, \dots, $G_m$ are pairwise disjoint, all of the crossings of $G$ contained in $G_1$, \dots, $G_m$.  

Suppose that $G$ isotopic to a spatial graph diagram $G'$.  Then each $G_i$ is isotopic to some $G_i'$ in $G'$, and $G_i'$ has at least as many crossings as $G_i$.  In addition, $G'$ may have crossings which are not contained in one of the $G_i'$.  Thus $G'$ has at least as many crossings as $G$.  \end{proof}

To see that we can apply this theorem to the spatial graph diagram in Figure~\ref{Framework}, observe that the orange and red links are reduced and alternating.  The blue is the Kinoshita $\theta$-curve which is known to have minimal crossing number $5$ \cite{7crossings}.  Finally, the purple is an embedding of the complete graph $K_6$ which was shown to have minimal crossing number in \cite{SimonWu}.  Thus $G$ has minimal crossing number by Theorem~\ref{theorem:framework}.

The following corollary is immediate by using Thistlethwaite's Theorem \cite{This}

\begin{corollary} \label{Corollary:framework} Let $G$ be a spatial graph diagram with a planar graph framework $P$, such that all of replacements for edges and vertices are adequate link diagrams.  Then $G$ has minimal crossing number.
\end{corollary}

\medskip

Let $G$ be a spatial graph diagram like the one in the Corollary which has even valence.  While there is a smoothing that gives us a pairwise disjoint collection of adequate links, not every smoothing necessarily produces an adequate link.  Thus $G$ is not necessarily adequate.

\begin{Open} Let $G$ be an adequate diagram of an even valent spatial graph.  Does $G$ necessarily have minimal crossing number?
\end{Open}

\medskip

\section{Rigid vertex spatial graphs}\label{rigid}
\begin{figure}[h!]
\includegraphics[scale=.3]{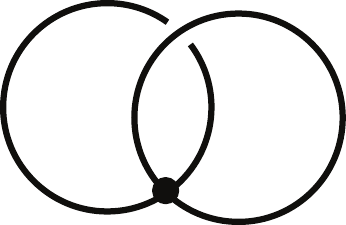}
\caption{This rigid vertex spatial graph has minimal crossing number $1$.}
\label{fig:onecrossing}
\end{figure}

 A \emph{rigid vertex spatial graph} is one which has a small planar disk around each vertex.  A \emph{rigid vertex isotopy} is an isotopy which preserves this set of disks.  Such an isotopy can turn over a disk, but the edges coming out of a disk cannot be permuted.  As a result, a diagram can have minimal crossing number as a rigid vertex graph, but not as a spatial graph.  For example Figure~\ref{fig:onecrossing} shows a spatial graph diagram that has minimal crossing number $1$ as a rigid vertex spatial graph, but would have minimal crossing number $0$ as an ordinary spatial graph.


Adams et al. \cite{Adams} defined reduced and alternating spatial graph diagrams as follows.

\begin{defi}\label{RedAltAdams} \cite{Adams} Let $G$ be a spatial graph diagram, which may or may not be a rigid vertex spatial graph.  
\begin{enumerate}
\item $G$ is {\bf reduced} if there is no circle $C$ in the plane of the diagram such that $C$ intersects $G$ transversely in a single crossing. 
\item $G$ is {\bf alternating} if 
\begin{itemize}
\item crossings alternate along each edge.
\item given be a planar neighborhood $N$ of vertices and uncrossed edges, first crossings along edges exiting any component of $N$ alternate as we go around the boundary of the component.
\end{itemize}
 \end{enumerate} 
\end{defi}

 \begin{AdamsT}\cite{Adams} Let $\Gamma$ be an even-valent rigid vertex spatial graph whose diagram $G$ is reduced and alternating.  If $G$ contains no uncrossed cycle, then $G$ has minimal crossing number.    \end{AdamsT}

\begin{figure}[h!]
			\centering	\includegraphics[width=.7\textwidth]{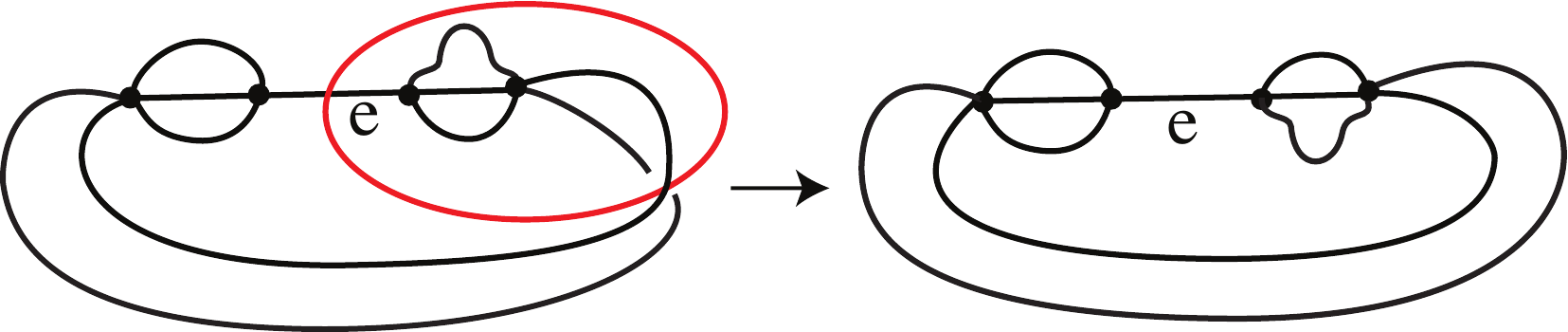}
			\caption{A rigid isotopy eliminates the crossing intersected by the red circle.}
			\label{Adams}
			\end{figure}

 Adams et al. shows that the above theorem for rigid-vertex graphs cannot be extended to graphs with at least one vertex of valence $3$, but asks whether it could be extended to graphs with at least one odd-valence vertex but none of valence $3$.  In particular, they ask for such an example which does not have a circle passing transversely through a single crossing and a single edge.  The spatial graph diagram in Figure~\ref{Adams} has such a circle (illustrated in red).  If there is such a circle going through a crossing and an edge $e$, then we can eliminate the crossing by a rigid isotopy that twists the graph around $e$ (see Figure~\ref{Adams}).  
  
\begin{figure}[h!]
			\centering	\includegraphics[width=.9\textwidth]{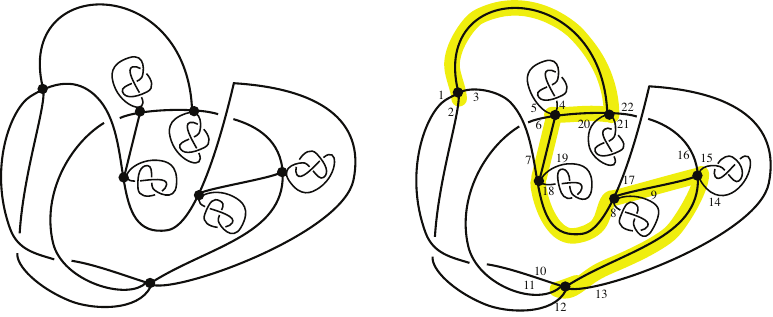}
			\caption{A spatial graph diagram which answers the question of Adams et al..  The neighborhood of the vertices and uncrossed edges is illustrated in yellow on the right.}
			\label{Definition}
			\end{figure}

 The diagram on the left in Figure~\ref{Definition} provides a rigid vertex spatial graph diagram which answers the question raised by Adams et al..  In particular, the neighborhood $N$ of vertices and uncrossed edges is illustrated in yellow on the right, with edges exiting $N$ numbered consecutively so that it is easy to check that first crossings alternate as we go along $\partial N$.  Observe that there are no vertices of valence $3$, no uncrossed cycles, and  no  circle passing transversely through a single crossing and a single edge.

\begin{figure}[h!]
			\centering	\includegraphics[width=\textwidth]{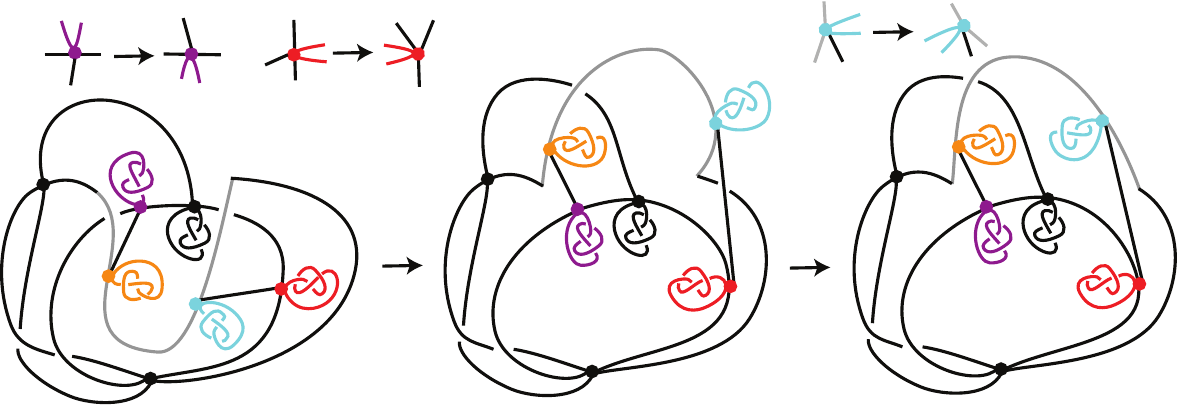}
			\caption{A rigid-vertex isotopy which turns over the red, purple, and turquoise disks reduces the number of crossings.}
			\label{New4_2}
			\end{figure}

Figure~\ref{New4_2} illustrates a rigid vertex isotopy which transforms the reduced alternating diagram on the left to a non-alternating diagram with fewer crossings.  In the first step we lift the grey arc, turning over the red and purple disks as indicated above the figure.  In the second step we turn over the turquoise disk and eliminate a crossing.  Note that while the initial diagram is alternating, neither the second nor the third diagram in Figure~\ref{New4_2} is alternating.  In particular, two consecutive edges exiting a neighborhood of the orange disk are overcrossings.

\bibliographystyle{amsplain}

  \bibliography{Splitting.bib}

@article {SimonWu,
    AUTHOR = {Flapan, Erica and Fletcher, Will and Nikkuni, Ryo},
     TITLE = {Reduced {W}u and generalized {S}imon invariants for spatial
              graphs},
   JOURNAL = {Math. Proc. Cambridge Philos. Soc.},
  FJOURNAL = {Mathematical Proceedings of the Cambridge Philosophical
              Society},
    VOLUME = {156},
      YEAR = {2014},
    NUMBER = {3},
     PAGES = {521--544},
      ISSN = {0305-0041,1469-8064},
   MRCLASS = {05C62 (05C10 57M15 57M25)},
  MRNUMBER = {3181638},
MRREVIEWER = {Jeremy\ L.\ Martin},
       DOI = {10.1017/S0305004114000073},
       URL = {https://doi.org/10.1017/S0305004114000073},
}

@article {7crossings,
   AUTHOR = {Moriuchi, Hiromasa},
     TITLE = {An enumeration of theta-curves with up to seven crossings},
   JOURNAL = {J. Knot Theory Ramifications},
  FJOURNAL = {Journal of Knot Theory and its Ramifications},
    VOLUME = {18},
      YEAR = {2009},
    NUMBER = {2},
     PAGES = {167--197},
      ISSN = {0218-2165,1793-6527},
   MRCLASS = {57M15 (05C10 57M25)},
  MRNUMBER = {2507922},
MRREVIEWER = {Chichen\ M.\ Tsau},
       DOI = {10.1142/S0218216509006884},
       URL = {https://doi.org/10.1142/S0218216509006884},
}

@article{Adams,
	author = {Adams, Colin and Dorman, Ryan and Foley, Kerryann and Kravis, Jonathan and Payne, Sam},
	doi = {10.1006/jctb.1999.1915},
	fjournal = {Journal of Combinatorial Theory. Series B},
	issn = {0095-8956},
	journal = {J. Combin. Theory Ser. B},
	mrclass = {05C10 (57M15)},
	mrnumber = {1710534},
	mrreviewer = {Lorenzo Traldi},
	number = {1},
	pages = {96--120},
	title = {Alternating graphs},
	url = {https://doi-org.ccl.idm.oclc.org/10.1006/jctb.1999.1915},
	volume = {77},
	year = {1999}}

@article {D,
    AUTHOR = {Diao, Yuanan},
     TITLE = {The additivity of crossing numbers},
   JOURNAL = {J. Knot Theory Ramifications},
  FJOURNAL = {Journal of Knot Theory and its Ramifications},
    VOLUME = {13},
      YEAR = {2004},
    NUMBER = {7},
     PAGES = {857--866},
      ISSN = {0218-2165,1793-6527},
   MRCLASS = {57M25},
  MRNUMBER = {2101230},
MRREVIEWER = {Masayuki\ Yamasaki},
       DOI = {10.1142/S0218216504003524},
       URL = {https://doi.org/10.1142/S0218216504003524},
}

@article {HY,
    AUTHOR = {Hara, Masao and Yamamoto, Makoto},
     TITLE = {Some links with nonadequate minimal-crossing diagrams},
   JOURNAL = {Math. Proc. Cambridge Philos. Soc.},
  FJOURNAL = {Mathematical Proceedings of the Cambridge Philosophical
              Society},
    VOLUME = {111},
      YEAR = {1992},
    NUMBER = {2},
     PAGES = {283--289},
      ISSN = {0305-0041,1469-8064},
   MRCLASS = {57M25},
  MRNUMBER = {1142747},
MRREVIEWER = {Mark\ E.\ Kidwell},
       DOI = {10.1017/S030500410007537X},
       URL = {https://doi.org/10.1017/S030500410007537X},
}

@article {KM,
    AUTHOR = {Kalfagianni, Efstratia and Mcconkey, Rob},
     TITLE = {Crossing numbers of cable knots},
   JOURNAL = {Bull. Lond. Math. Soc.},
  FJOURNAL = {Bulletin of the London Mathematical Society},
    VOLUME = {56},
      YEAR = {2024},
    NUMBER = {11},
     PAGES = {3400--3411},
      ISSN = {0024-6093,1469-2120},
   MRCLASS = {57K10 (57K14 57K16)},
  MRNUMBER = {4828022},
       DOI = {10.1112/blms.13140},
       URL = {https://doi.org/10.1112/blms.13140},
}

@article{Kauff,
	author = {Kauffman, Louis H.},
	doi = {10.1016/0040-9383(87)90009-7},
	fjournal = {Topology. An International Journal of Mathematics},
	issn = {0040-9383},
	journal = {Topology},
	mrclass = {57M25},
	mrnumber = {899057},
	mrreviewer = {J. S. Birman},
	number = {3},
	pages = {395--407},
	title = {State models and the {J}ones polynomial},
	url = {https://doi-org.ccl.idm.oclc.org/10.1016/0040-9383(87)90009-7},
	volume = {26},
	year = {1987}}

@article {Lack,
    AUTHOR = {Lackenby, Marc},
     TITLE = {The crossing number of composite knots},
   JOURNAL = {J. Topol.},
  FJOURNAL = {Journal of Topology},
    VOLUME = {2},
      YEAR = {2009},
    NUMBER = {4},
     PAGES = {747--768},
      ISSN = {1753-8416,1753-8424},
   MRCLASS = {57M25},
  MRNUMBER = {2574742},
       DOI = {10.1112/jtopol/jtp028},
       URL = {https://doi.org/10.1112/jtopol/jtp028},
}

@article{LT,
	author = {Lickorish, W. B. R. and Thistlethwaite, M. B.},
	doi = {10.1007/BF02566777},
	fjournal = {Commentarii Mathematici Helvetici},
	issn = {0010-2571},
	journal = {Comment. Math. Helv.},
	mrclass = {57M25},
	mrnumber = {966948},
	mrreviewer = {G. Burde},
	number = {4},
	pages = {527--539},
	title = {Some links with nontrivial polynomials and their crossing-numbers},
	url = {https://doi-org.ccl.idm.oclc.org/10.1007/BF02566777},
	volume = {63},
	year = {1988}}

@article{Mura1,
	author = {Murasugi, Kunio},
	doi = {10.1016/0040-9383(87)90058-9},
	fjournal = {Topology. An International Journal of Mathematics},
	issn = {0040-9383},
	journal = {Topology},
	mrclass = {57M25},
	mrnumber = {895570},
	mrreviewer = {Louis H. Kauffman},
	number = {2},
	pages = {187--194},
	title = {Jones polynomials and classical conjectures in knot theory},
	url = {https://doi-org.ccl.idm.oclc.org/10.1016/0040-9383(87)90058-9},
	volume = {26},
	year = {1987}}

@article {Mura,
    AUTHOR = {Murasugi, Kunio},
     TITLE = {On the braid index of alternating links},
   JOURNAL = {Trans. Amer. Math. Soc.},
  FJOURNAL = {Transactions of the American Mathematical Society},
    VOLUME = {326},
      YEAR = {1991},
    NUMBER = {1},
     PAGES = {237--260},
      ISSN = {0002-9947,1088-6850},
   MRCLASS = {57M25},
  MRNUMBER = {1000333},
MRREVIEWER = {Hugh\ Reynolds\ Morton},
       DOI = {10.2307/2001863},
       URL = {https://doi.org/10.2307/2001863},
}

@article{Sawollek2,
	author = {Sawollek, J\"{o}rg},
	doi = {10.1016/S0166-8641(97)00268-X},
	fjournal = {Topology and its Applications},
	issn = {0166-8641},
	journal = {Topology Appl.},
	mrclass = {57M15 (05C10)},
	mrnumber = {1688477},
	number = {3},
	pages = {261--273},
	title = {Alternating diagrams of {$4$}-regular graphs in {$3$}-space},
	url = {https://doi-org.ccl.idm.oclc.org/10.1016/S0166-8641(97)00268-X},
	volume = {93},
	year = {1999}}

@article{Stoimenow1,
	author = {Stoimenow, Alexander},
	doi = {10.1515/forum-2011-0121},
	fjournal = {Forum Mathematicum},
	issn = {0933-7741},
	journal = {Forum Math.},
	mrclass = {57M25 (57M27)},
	mrnumber = {3228928},
	mrreviewer = {Kyoung Il Park},
	number = {4},
	pages = {1187--1246},
	title = {On the crossing number of semiadequate links},
	url = {https://doi-org.ccl.idm.oclc.org/10.1515/forum-2011-0121},
	volume = {26},
	year = {2014}}

@article{Thistle,
	author = {Thistlethwaite, Morwen B.},
	doi = {10.1016/0040-9383(87)90003-6},
	fjournal = {Topology. An International Journal of Mathematics},
	issn = {0040-9383},
	journal = {Topology},
	mrclass = {57M25},
	mrnumber = {899051},
	mrreviewer = {G. Burde},
	number = {3},
	pages = {297--309},
	title = {A spanning tree expansion of the {J}ones polynomial},
	url = {https://doi-org.ccl.idm.oclc.org/10.1016/0040-9383(87)90003-6},
	volume = {26},
	year = {1987}}

@article{This,
	author = {Thistlethwaite, Morwen B.},
	doi = {10.1007/BF01394334},
	fjournal = {Inventiones Mathematicae},
	issn = {0020-9910},
	journal = {Invent. Math.},
	mrclass = {57M25},
	mrnumber = {948102},
	mrreviewer = {V. G. Turaev},
	number = {2},
	pages = {285--296},
	title = {On the {K}auffman polynomial of an adequate link},
	url = {https://doi-org.ccl.idm.oclc.org/10.1007/BF01394334},
	volume = {93},
	year = {1988}}

\end{document}